\documentclass{amsart}
\usepackage{subfigure}
\usepackage{amsmath}
\usepackage{amsfonts}
\usepackage{graphics}
\usepackage{pstricks}

\title{Where are the Natural Numbers in Hilbert's Foundations of Geometry?}
\author{Phil Scott and Jacques D Fleuriot}

\address{Artificial Intelligence and its Applications Institute (AIAI)\\ 
            School of Informatics\\
            University of Edinburgh\\
            10 Crichton Street \\
            EH8 9AB \\
            United Kingdom
            }
            
\email{jacques.fleuriot@ed.ac.uk}

\begin{document}
\maketitle

\begin{abstract}
Hilbert's \emph{Foundations of Geometry} was perhaps one of the most influential works of geometry in the 20th century and its axiomatics was the first systematic attempt to clear up the logical gaps of the \emph{Elements}. But does it have gaps of its own? In this paper, we discuss a logical issue, asking how Hilbert is able to talk about natural numbers within a foundational synthetic geometry. We clarify the matter, showing how to obtain the natural numbers using a very modest subset of his axioms.
\end{abstract}

\section{Introduction}
It has been argued that the Pythagoreans took their philosophical foundations to be the natural numbers~\cite{EvolutionEuclideanElements}, and if so, the discovery of irrational numbers marks the first crisis in the foundations of mathematics. By the time of Euclid, the theory of natural numbers was grounded instead on supposedly more secure geometric notions, and placed on par with an equally rich theory of rational and irrational magnitudes.

We see these ideas in Euclid's \emph{Elements}~\cite{HeathElements}, possibly one of the most influential books in history~\cite{BoyerEuclidInfluence}, and possibly the most famous textbook in mathematics. Its presentation of axiomatic geometry went largely unquestioned until the 19th century. Then, as mathematicians were becoming more attentive to logical principles, Pasch and Hilbert concluded that the text contained missing assumptions. To be truly rigorous, Pasch had declared that the process of inference in all proofs should be entirely independent of the meaning of any geometrical term~\cite{ModernAxiomatics}, while David Hilbert, anticipating his own ``formalist programme'', is famously reported to have said that ``mug'', ``table'', and ``chair'' should be substitutes for ``point'', ``line'' and ``plane''~\cite{TableChairMug}. The conclusion of this view was Hilbert's \emph{Foundations of Geometry}, in which Euclid's five axioms became nineteen axioms, organised into five groups.

As Poincar\'{e} explained in his review of the first edition of the \emph{Foundations of Geometry}~\cite{PoincareReview}, we can understand this idea of rigour in terms of a purely mechanical symbolic machine. If one feeds the axioms in as input, then all of Hilbert's theorems should be delivered as output. Thus, we guarantee that intuition does not hide implicit axioms relating the primitive concepts, since a purely mechanical machine has no such intuitions.

Such a mechanical symbolic approach to mathematics was made in earnest by Russell and Whitehead when they wrote their \emph{Principia}, albeit without the aid of a real logic machine. The task was clearly laborious. Poincar\'{e} remarked derisively that ``if it requires 27 equations to establish that 1 is a number, how many will it require to demonstrate a real theorem?''~\cite{PoincareShackles}. Russell admitted that he never fully recovered from the effort. 

Now with modern computers and automated theorem provers, the effort is greatly alleviated. Machines can now check enormous symbolic proofs far more quickly and reliably than any human. Furthermore, the machines assist us with book-keeping and cross-referencing, and relieve us the burden of generating the tedious and simple details that Poincar\'{e} presumably considered  ``puerile''. In fact, with the invention of Herbert Simon's mechanical \emph{Logic Theorist}, Russell was forced to reflect on ten years writing proofs by hand as a ``wasted'' effort~\cite{SimonObituary}. For two excellent introductions to computer assisted formalisation, suitable for non-experts, see Hales~\cite{HalesFormalProof} and Harrison~\cite{FormalizedMathematics}.

In keeping with Poincar\'{e}'s mechanistic conclusion of logical rigour, we now assure the reader that the work in this paper has been backed up by symbolic definitions and proofs that have been fully formalised and machine-checked using the theorem prover HOL Light \cite{HOLLight}. However, for readability, we present our results in the more familiar, and perhaps less ``puerile'', mathematical vernacular.

Now as we found when we attempted to formalise Hilbert's theory, there are some subtle representation issues and choices to be made concerning the underlying logic~\cite{MeikleFleuriotFormalizingHilbert}. In this paper, we are concerned with just one issue: the logical status of the natural numbers.

\section{$\mathbb{N}$ in \emph{The Foundations of Geometry}}
Definition 4 of Book V of the \emph{Elements} reads:

\begin{quote}Magnitudes are said to \emph{have a ratio} to one another which can, when multiplied, exceed one another.
\end{quote}

The way this is used makes it clear that it should have been an axiom, an issue that was corrected in the \emph{Foundations of Geometry}, where it appears as Axiom~V,~1 along with an axiom of completeness.

\begin{quote}
\noindent V,~1 (Axiom of measure or Archimedes' Axiom). If $AB$ and $CD$ are any segments then there exists a number $n$ such that $n$ segments $CD$ constructed contiguously from $A$, along the ray from $A$ through $B$, will pass beyond the point $B$.

\noindent V,~2 (Axiom of line completeness). An extension of a set of points on a line with its order and congruence relations that would preserve the relations existing among the original elements as well as the fundamental properties of line order and congruence that follows from Axioms I-III, and from V,~1 is impossible.
\end{quote}

In this last group of axioms in Hilbert's system, the assumed logic is substantially stronger than that assumed in the earlier groups. In the first axiom, we are told there exist natural numbers dependent on segments. Natural numbers have not been defined at this point in the text, so here, they are being treated as a logical primitive. This matter seems incongruent with Euclid's \emph{Elements}, which treated the natural numbers geometrically. 

And in Axiom~V,~2, Hilbert is quantifying over sets of points. But the combination of a second-order logic with a theory of natural numbers gives us second-order arithmetic, wherein we can obtain a model of the \emph{real numbers}. What is the \emph{logical priority} here? Is the basic logic here as strong as second-order arithmetic, or should we try to recover natural numbers from geometrical notions? To answer this question, we note that Hilbert wanted to avoid the complexity of analysis and produce a strictly synthetic theory~\cite{MajerHilbertKleinComplexity}, and he went to great length in later chapters to define an arithmetic based on magnitudes. Thus, it would seem we should define the natural numbers geometrically. The question then is \emph{where} would the definition appear?

\section{Group~I}
The first and largest group contains eight axioms to define properties of \emph{incidence}. There are three primitive domains, namely \emph{points}, \emph{lines} and \emph{planes}, and two primitive relations, one for point-line incidence, and another for point-plane incidence. 

\begin{figure}
\begin{minipage}[t]{3cm}
\scalebox{0.3}{\includegraphics{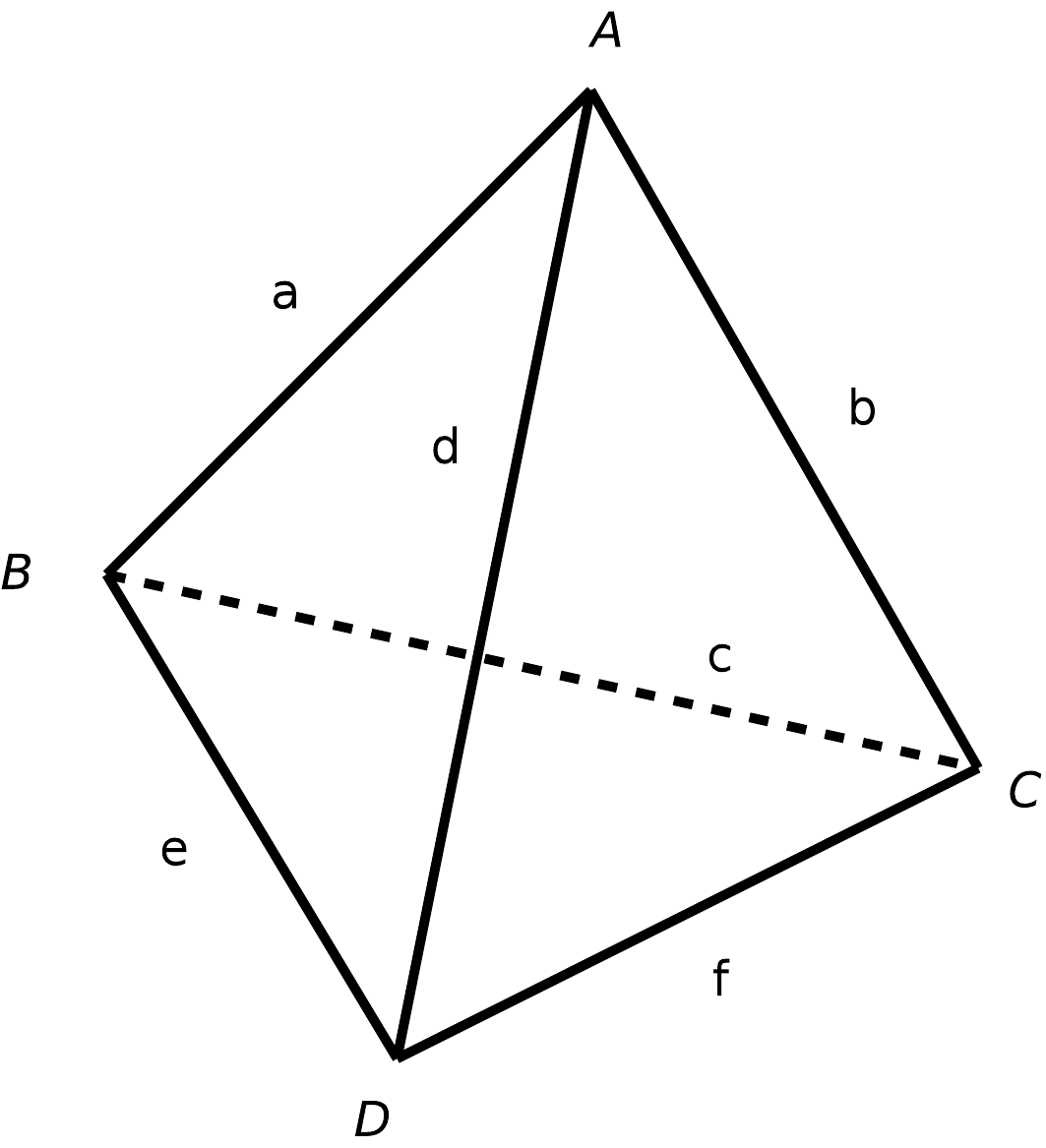}}
\end{minipage}
\begin{minipage}[b]{6cm}
\begin{itemize}
\item $A$ and $B$ lie on $a$;
\item $A$ and $C$ lie on $b$;
\item $B$ and $C$ lie on $c$;
\item $A$ and $D$ lie on $d$;
\item $B$ and $D$ lie on $e$;
\item $C$ and $D$ lie on $f$.
\end{itemize}
\end{minipage}
\\
\begin{itemize}
\item $A$, $B$ and $C$ lie on plane $\alpha$;
\item $A$, $B$ and $D$ lie on plane $\beta$;
\item $A$, $C$ and $D$ lie on plane $\gamma$;
\item $B$, $C$ and $D$ lie on plane $\delta$.
\end{itemize}

\caption{Finite Model of Group~I}\label{fig:SmallestModel}
\end{figure}

It turns out that the axioms governing these two relations can be minimally realised in four vertices $A$, $B$, $C$, $D$, six lines $a$, $b$, $c$, $d$, $e$, $f$, and the four planes $ABC$, $ABD$, $ACD$ and $BCD$ of a tetrahedron (see Figure~\ref{fig:SmallestModel}). That is, we can prove that 14 objects in certain relations satisfy the axioms, and prove that the axioms commit us to at least these 14 objects.

Since finite models are possible for the first group of axioms, we know we cannot derive the existence of an infinite collection such as the natural numbers. Thus, we look to the second group of axioms.

\section{Group~II}
The next group of axioms concern \emph{order} and \emph{orientation} in terms of a primitive three-place predicate of betweenness, with which we can state that one point lies strictly between two others along a line. The axioms to govern this notion are given in Figure \ref{fig:GroupII}.

\begin{figure}
  \begin{enumerate}
  \item[II,1] If a point $B$ lies between a point $A$ and a point $C$ then the points $A$, $B$, $C$ are three distinct points of a line, and $B$ then also lies between $C$ and $A$.
  \item[II,2] For two points $A$ and $C$, there always exists at least one point $B$ on the line $AC$ such that $C$ lies between $A$ and $B$.
  \item[II,3] Of any three points on a line there exists no more than one that lies between the other two.
  \item[II,4] Let $A$, $B$, $C$ be three points that do not lie on a line and let $a$ be a line in the plane $ABC$ which does not meet any of the points $A$, $B$, $C$. If the line $a$ passes through a point of the segment $AB$, it also passes through a point of the segment $AC$, or through a point of the segment $BC$.
  \end{enumerate}
  \caption{Group II axioms}
  \label{fig:GroupII}
\end{figure}

These axioms have been heavily revised since the first edition as various redundancies were discovered. Originally, Axiom~II,~2 had an additional clause asserting that we can always find a point between two others. This turned out to be derivable and is Theorem~3 in the tenth edition. Axiom~II,~2 is then effectively a rigorous formulation of the second postulate of the \emph{Elements}, allowing us to continue a finite straight line (or line segment) along a straight line. Unlike Euclid's formulation, the use of the betweenness relation makes it clear that this is an axiom of \emph{order}.

In the first edition, there were further redudancies and these too were later removed\footnote{This despite Hilbert's claim in the metatheory of his first edition that his axioms were clearly independent!}. The quantification in Axiom~II,~3 was stronger, claiming that, given three points, there is \emph{exactly} one that lies between the other two. By the tenth edition, the existence of the point was derived as Theorem~4 in a proof attributed to Wald~\cite{FoundationsOfGeometry}. Finally, a ``transitivity'' axiom was removed to become Theorem~5, the proof due to Moore~\cite{MooreProof}. The theorems are given in Figure~\ref{fig:Theorems345}.

\begin{figure}
  \begin{quote}
    \noindent THEOREM 3. For any two points $A$ and $C$ there always exists at least one point $D$ on the line $AC$ that lies between $A$ and $C$.
    
    \noindent THEOREM 4. Of any three points $A$, $B$, $C$ on a line there always is one that lies between the other two.
    
    \noindent THEOREM 5. Given any four points on a line, it is always possible to label them $A$, $B$, $C$, $D$ in such a way that the point labeled $B$ lies between $A$ and $C$ and also between $A$ and $D$, and furthermore, that the point labeled $C$ lies between $A$ and $D$ and also between $B$ and $D$.

    \noindent THEOREM 6 (generalization of Theorem 5). Given any finite number of points on a line it is always possible to label them $A$, $B$, $C$, $D$, $E$, $\ldots$, $K$ in such a way that the point labelled $B$ lies between $A$ and $C$, $D$, $E$, $\ldots$, $K$, the point labelled $C$ lies between $A$, $B$, and $D$, $E$, $\ldots$, $K$, $D$ lies between $A$, $B$, $C$ and $E$, $\ldots$, $K$, etc. Besides this order of labelling there is only the reverse one that has the same property.

    \noindent THEOREM 7. Between any two points on a line there exists an infinite number of points. 
  \end{quote}
\caption{Some Group II Theorems}
\label{fig:Theorems345}
\end{figure}

\begin{figure}

\scalebox{0.3}{\includegraphics{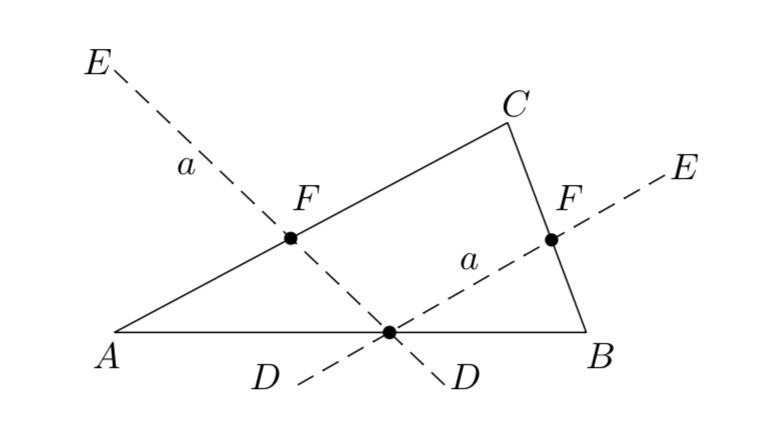}}
\caption{Axiom~II,~4}
\label{fig:PaschDiagram}
\end{figure}

After revising, the Group~II axioms consist of constraints on the between relation, an axiom to extend a line segment in a given direction, and finally a planar axiom of order. This last axiom, depicted in Figure~\ref{fig:PaschDiagram}, was identified by Pasch as a crucial missing postulate from Euclid's~Elements. It is used heavily in the proofs of Theorem~4 and Theorem~5, each of which starts with a linear assumption, builds a planar diagram, and then uses Axiom~II,~4 repeatedly to produce a linear conclusion.

\section{Infinity}\label{sec:Infinity}
Hilbert does not give a proof for Theorem~7, but it claims here that the axioms now only have infinite models. Indeed, we could presumably apply Theorem~3 repeatedly to obtain as many points as desired between two others. We start with distinct points $A$ and $C$, and then obtain points $D$, $E$, $F$, $\ldots$, $Y$, $Z$ such that the point $D$ is between $A$ and $C$, the point $E$ is between $A$ and $D$, the point $F$ is between $A$ and $E,\ldots$ and the point $Y$ is between $A$ and $Z$.

We just need to prove that these points are distinct, and for this, we can use  Theorem~5. We know that if $D$ is between $A$ and $C$ and $E$ is between $A$ and $D$, that $E$ must also be between $A$ and $C$. And thus, we know that these points are distinct. Repeating this line of reasoning therefore gives us a \emph{potentially} infinite number of distinct points.

This is not enough to give us the natural numbers unless we can logically ``carve out'' the results of this process into something which we can quantify over. For us, it suffices that we can exhibit such a process as a function. 

\section{A Geometric Successor}
\begin{figure}
\scalebox{0.25}{\includegraphics{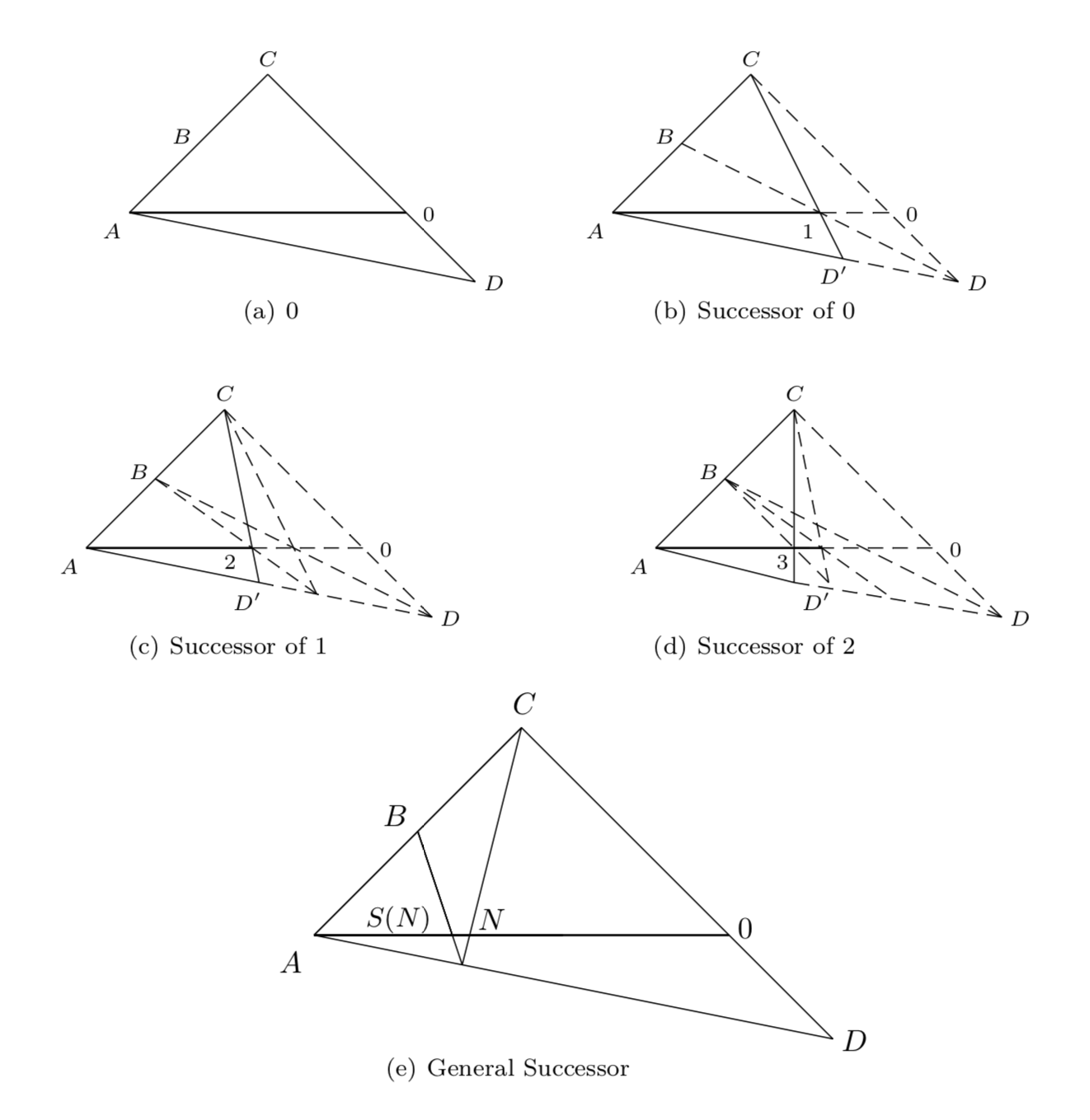}}
\caption{The successor function}
\label{fig:successor}
\end{figure}

The function we need is the traditional successor function, and it is possible to witness such a function as a unique object in Group~II. We have chosen our witness to have two properties: firstly, it can be understood as a formalisation of Hilbert's Theorem~7; secondly, the points obtained by the successor function are uniquely defined (unlike the potential infinity of points from the previous section, which we would have to obtain via the axiom of choice). The witness is based on the diagram in Figure~\ref{fig:successor}.\footnote{Note that in the diagram, we could have switched the labels, using $0$, $1$, $2$, $\ldots$ to denote the moving point $D'$.}

Formally, each diagram is a set of points $A$, $B$, $C$, $D$, $N$ and $0$, satisfying the constraints in Figure~\ref{fig:DiagramConstraints}. The fact that points exist satisfying these constraints is proven in Theorem~3, which uses a diagram exactly as shown in Figure~\ref{fig:successor}(a).

\begin{figure}
  \begin{enumerate}
  \item the points $A$, $B$ and $0$ are not collinear; 
  \item the point $B$ lies between $A$ and $C$;
  \item the point $0$ lies between $C$ and $D$;
  \item either:
    \begin{enumerate}
    \item the point $N$ lies between $A$ and $0$, or
    \item $N = 0$
    \end{enumerate}
  \end{enumerate}
\label{fig:DiagramConstraints}
\end{figure}

Each diagram is represented by six points. The points $A$, $B$, $C$, $D$ and $0$ are fixed under the successor function, while the point $N$ is moved closer to $A$. Note, however, that at this stage, we only have axioms governing incidence and orientation. In particular, the first two groups of axioms are consistent with a geometry allowing infinitely small distances. So, if the point $B$ was infinitely close to $C$, then the points $1,2,3,\ldots$ would be infinitely close to $0$. This does not affect the fact that we have a successor function, but is worth bearing in mind nevertheless: when it comes to working in the foundations of geometry, diagrams are often misleading. They represent many more assumptions than are normally in force.

The set of these diagrams is not, of course, provably isomorphic to the natural numbers. In the model of Euclidean geometry, for instance, the set of geometrically \emph{similar}\footnote{Note that we do not have enough axioms at this stage to define this notion.} diagrams is isomorphic to an interval of \emph{real numbers}. To obtain the naturals, we use induction as in set theory or simple type theory~\cite{ChurchTheoryOfTypes}, whereby we intersect all sets of diagrams which contain a chosen $0$ and are closed under our successor function. It is then sufficient to prove that this successor function is one-one and not onto (something normally taken as axiomatic).

Informally, and for the purposes of explaining our formalisation, we will identify diagrams by the point $N$. In this way, when we talk of the object $0$, we may be referring to the \emph{diagram} 0, which are the six points shown in Figure~\ref{fig:successor}(a) or to the \emph{point} 0, which is just one of these six points. Similarly, we shall talk about the \emph{diagram} that is the successor of 0, as well as the \emph{point} that is the successor of 0.\footnote{In our machine-checked formalisation, these matters are necessarily unambiguous.}

Thus, the successor of $0$ is the diagram obtained by replacing $0$ with $1$, the intersection of $BD$ and $A0$. To obtain the next successor, we first find the intersection of $C1$ and $AD$, namely the point $D'$. We then replace $1$ with the intersection of $BD'$ and $A0$. 

In general, the successor of a diagram is obtained by finding $D'$, the intersection of $CN$ and $AD$, and then finding the intersection of $BD'$ and $A0$. Formally, if $A$, $B$, $C$, $D$, $0$ and $N$ are a diagram, then the successor to this diagram are the points $A$, $B$, $C$, $D$, $0$ and that provably unique $S$ such that there is a point $D'$ with:

\begin{enumerate}
\item the points $C$, $D'$ and $N$ collinear;
\item the points $B$, $D'$ and $S$ collinear;
\item the points $A$, $D$ and $D'$ collinear;
\item $S$ lying between $A$ and $0$;
\item $S$ lying between $A$ and $N$.
\end{enumerate}

We just need to show that these points satisfy the constraints given in Figure~\ref{fig:DiagramConstraints}. The key step is to notice that the diagrams involving the points $A,B,C,0,D$ and $A,B,C,N,D'$ satisfy the same constraints. In each case, we apply Axiom~II,~4 to $\triangle AC0$ and the line $BD'$ to locate a point $S$ between $A$ and $0$. For our initial diagram, we can apply this argument directly. For the other diagrams, we just need to find the point $D'$.

To do this, we apply Axiom~II,~4 to $\triangle A0D$ and the line $BN$, to place the point $D'$ between $A$ and $D$. We can now use the previous argument to locate $S$ between $A$ and $N$. Finally, since $N$ is between $A$ and $0$, Theorem~5 shows that the point $S$ must also lie between $A$ and $0$.

With our function defined, we must show that it satisfies the crucial property of being a successor function: namely, that it is one-one but not onto (and thus that the set of points in its image are Dedekind-infinite). We can show that the function is one-one by reasoning about incidence alone, and it is clearly not onto since the point $0$ defines the first diagram, while all images of our successor function use a point $S$ which is strictly between $A$ and $0$.

\section{Conclusion}
We have derived the existence of a set of natural numbers using just the first two groups of Hilbert's axiom system. Ordinarily, the existence of such a set would be taken as axiomatic, as is the case in ZF set theory\footnote{Strictly speaking, the axiom of infinity in ZF set theory merely says there is a set containing zero and closed under successor.}. But this is not necessary, and was sometimes considered logically questionable. In the \emph{Principia}, Russell declared that the existence of an actual infinity of natural numbers was an empirical, not logical matter. Hilbert himself would later argue that the infinite is merely a useful device, and not a logical given~\cite{OnInfinite}.

But in Hilbert's geometry, we can understand natural numbers as just another geometric concept. They are, in our formalism, a set of diagrams involving a triangle, each generated from its predecessor by applying Axiom~II,~4. Thus, when Hilbert mentions natural numbers in his text (such as when formalising Axiom~V,~1) we can understand him to be implicitly using a derivative geometric notion, one that we have now explicated based on his first two groups of axioms.

This makes explicit  a gap in the \emph{Foundations of Geometry}, and helps us realise synthetic geometry as a true foundation for number theory, as it was understood by its predecessor, the \emph{Elements}. 


\bibliographystyle{plain}
\bibliography{infinity}
\end{document}